\documentclass[12pt,a4paper]{article}
\usepackage{epsfig,latexsym,amsfonts,amssymb,amsmath,amscd,
epic,graphics,theorem}%,pxfonts}
\theoremstyle{plain}
\newtheorem{lemma}{Lemma}

\newtheorem{remark}[lemma]{Remark}
\newtheorem{example}[lemma]{Example}
\newtheorem{theorem}[lemma]{Theorem}
\newtheorem{definition}[lemma]{Definition}
\newtheorem{corollary}[lemma]{Corollary}

\newtheorem{question}[lemma]{Question}
{\theorembodyfont{\rmfamily}  \font\ncsc=cmcsc10
 \font\ntt=cmtt12

\begin{document}
\newcommand{\pperp}{\hbox{$\perp\hskip-6pt\perp$}}
\newcommand{\ssim}{\hbox{$\hskip-2pt\sim$}}
\newcommand{\aleq}{{\ \stackrel{3}{\le}\ }}
\newcommand{\ageq}{{\ \stackrel{3}{\ge}\ }}
\newcommand{\aeq}{{\ \stackrel{3}{=}\ }}
\newcommand{\bleq}{{\ \stackrel{n}{\le}\ }}
\newcommand{\bgeq}{{\ \stackrel{n}{\ge}\ }}
\newcommand{\beq}{{\ \stackrel{n}{=}\ }}
\newcommand{\cleq}{{\ \stackrel{2}{\le}\ }}
\newcommand{\cgeq}{{\ \stackrel{2}{\ge}\ }}
\newcommand{\ceq}{{\ \stackrel{2}{=}\ }}
\newcommand{\N}{{\mathbb N}}
\newcommand{\A}{{\mathbb A}}
\newcommand{\K}{{\mathbb K}}
\newcommand{\Z}{{\mathbb Z}}
\newcommand{\R}{{\mathbb R}}
\newcommand{\C}{{\mathbb C}}
\newcommand{\Q}{{\mathbb Q}}
\newcommand{\PP}{{\mathbb P}}
\newcommand{\Pic}{{\operatorname{Pic}}}\newcommand{\Sym}{{\operatorname{Sym}}}
\newcommand{\oeps}{{\overline\eps}}\newcommand{\idim}{{\operatorname{idim}}}
\newcommand{\oDel}{{\widetilde\Del}}
\newcommand{\cD}{{\mathcal D}}\newcommand{\mt}{{\operatorname{mt}}}
\newcommand{\ord}{{\operatorname{ord}}}\newcommand{\Id}{{\operatorname{Id}}}
\newcommand{\Span}{{\operatorname{Span}}}\newcommand{\As}{{\operatorname{As}}}
\newcommand{\Ker}{{\operatorname{Ker}}}\newcommand{\nd}{{\operatorname{nd}}}
\newcommand{\Ann}{{\operatorname{Ann}}}
\newcommand{\Fix}{{\operatorname{Fix}}}
\newcommand{\sign}{{\operatorname{sign}}}
\newcommand{\Tors}{{\operatorname{Tors}}}
\newcommand{\Ima}{{\operatorname{Im}}}
\newcommand{\oi}{{\overline i}}
\newcommand{\oj}{{\overline j}}
\newcommand{\ob}{{\overline b}}
\newcommand{\os}{{\overline s}}
\newcommand{\oa}{{\overline a}}
\newcommand{\oy}{{\overline y}}
\newcommand{\ow}{{\overline w}}
\newcommand{\ot}{{\overline t}}
\newcommand{\oz}{{\overline z}}
\newcommand{\eps}{{\varepsilon}}
\newcommand{\proofend}{\hfill$\Box$\bigskip}
\newcommand{\Int}{{\operatorname{Int}}}
\newcommand{\pr}{{\operatorname{pr}}}
\newcommand{\Hom}{{\operatorname{Hom}}}
\newcommand{\rk}{{\operatorname{rk}}}\newcommand{\Ev}{{\operatorname{Ev}}}
\newcommand{\im}{{\operatorname{Im}}}
\newcommand{\sk}{{\operatorname{sk}}}\newcommand{\DP}{{\operatorname{DP}}}
\newcommand{\const}{{\operatorname{const}}}
\newcommand{\Sing}{{\operatorname{Sing}}\hskip0.06cm}
\newcommand{\conj}{{\operatorname{Conj}}}
\newcommand{\Cl}{{\operatorname{Cl}}}
\newcommand{\Crit}{{\operatorname{Crit}}}
\newcommand{\Ch}{{\operatorname{Ch}}}
\newcommand{\discr}{{\operatorname{discr}}}
\newcommand{\Tor}{{\operatorname{Tor}}}
\newcommand{\Conj}{{\operatorname{Conj}}}
\newcommand{\vol}{{\operatorname{vol}}}
\newcommand{\defect}{{\operatorname{def}}}
\newcommand{\codim}{{\operatorname{codim}}}
\newcommand{\ov}{{\overline v}}
\newcommand{\ox}{{\overline{x}}}
\newcommand{\bw}{{\boldsymbol w}}
\newcommand{\bx}{{\boldsymbol x}}
\newcommand{\bd}{{\boldsymbol d}}
\newcommand{\bz}{{\boldsymbol z}}\newcommand{\bp}{{\boldsymbol p}}
\newcommand{\tet}{{\theta}}
\newcommand{\Del}{{\Delta}}
\newcommand{\bet}{{\beta}}
\newcommand{\kap}{{\varkappa}}
\newcommand{\del}{{\delta}}
\newcommand{\sig}{{\sigma}}
\newcommand{\alp}{{\alpha}}
\newcommand{\Sig}{{\Sigma}}
\newcommand{\Gam}{{\Gamma}}
\newcommand{\gam}{{\gamma}}
\newcommand{\Lam}{{\Lambda}}
\newcommand{\lam}{{\lambda}}
\newcommand{\SC}{{SC}}
\newcommand{\MC}{{MC}}
\newcommand{\nek}{{,...,}}
\newcommand{\cim}{{c_{\mbox{\rm im}}}}
\newcommand{\clM}{\tilde{M}}
\newcommand{\clV}{\bar{V}}

%for mnotes
{\catcode`\@=11
\gdef\n@te#1#2{\leavevmode\vadjust{%
 {\setbox\z@\hbox to\z@{\strut#1}%
  \setbox\z@\hbox{\raise\dp\strutbox\box\z@}\ht\z@=\z@\dp\z@=\z@%
  #2\box\z@}}}
\gdef\leftnote#1{\n@te{\hss#1\quad}{}}
\gdef\rightnote#1{\n@te{\quad\kern-\leftskip#1\hss}{\moveright\hsize}}
\gdef\?{\FN@\qumark}
\gdef\qumark{\ifx\next"\DN@"##1"{\leftnote{\rm##1}}\else
 \DN@{\leftnote{\rm??}}\fi{\rm??}\next@}}
\def\mnote#1{\leftnote{\vbox{\hsize=2.5truecm\footnotesize
\noindent #1}}}

\title{Duality of planar and spacial curves: new insight}
\author{Victor Kulikov and Eugenii Shustin}
\date{}
\maketitle
\begin{abstract}
We study the geometry of equiclassical strata of the discriminant in the space of
plane curves of a given degree, which are families of curves of given degree, genus and class
(degree of the dual curve). Our main observation is that the use of duality
transformation leads to a series of new sufficient conditions for a regular
behavior of the equiclassical stratification.
We also discuss duality of curves in higher-dimensional projective spaces and in Grassmannians
with focus on similar questions of the regularity of equiclassical families of spacial curves.
\end{abstract}

{\small

}

\section*{Introduction}
Our main goal is to show that the issue of duality of planar and spacial projective curves
is not exhausted yet, and it still yields new results and raises interesting open problems.

The duality of plane algebraic curves is a classically known relation, which
has been used in the study of plane algebraic curves since 19th century. It is tightly linked to
the stratification of the discriminant in the space of curves of a given degree into equiclassical
strata, i.e. families of curves of given degree, genus and class (degree of the dual curve).
However, previous considerations of the geometry of equiclassical strata (see \cite{DH,Sh1})
have not taken the duality relation into account. We address regularity properties of equiclassical strata,
like that to be locally transverse intersection, to be irreducible, to have a nodal-cuspidal general member,
and we show that a combination of duality isomorphisms with tools of deformation theory leads to series
of new numerical sufficient conditions for the above regularity properties. The results of our study are
absorbed in Theorems \ref{tec3}, \ref{tlc1}, \ref{tlc2}, \ref{tqc1}, and \ref{tqc2} in Section \ref{sec2}.
For example, by \cite[Theorem 1.2]{DH} and \cite[Theorem 1.1]{Sh1}, the equiclassical family $V_{d,g,c}$
parameterizing irreducible curves of degree $d$, genus $g$,
and class $c$ has expected dimension $d-g+c+1$ and its generic member is a curve
with ordinary nodes and cusps if $\kap-2\delta\le3d-4$, where $\kap=d(d-1)-c$, $2\delta=(d-1)(d-2)-2g$. We show that
$V_{d,g,c}$ possesses the above properties as long as $5\kap-6\delta\le d^2+6d-3$
(cf. Theorem \ref{tqc1}(2) in Section \ref{sec2}).

The duality between non-degenerate curves in higher-dimensional spaces was discovered by Ragni Piene
\cite[\S 5]{P}, and later extended to curves in Grassmannians \cite{Per}. In Section \eqref{sec3},
we provide a brief account of this duality relation with elementary proofs, introduce
equiclassical families of spacial curves and formulate several open questions on
their geometry and on duality morphisms between them.

We work over the complex field, though most of results can be stated over any
algebraically closed field of characteristic zero.

{\small
{\bf Acknowledgements}.
The first author has been supported
by grants of  RFBR 14-01-00160 and 15-01-02158, and by the Government of the Russian Federation within the framework of the implementation of the 5-100 Programme Roadmap of the National Research University  Higher School of Economics, AG Laboratory. The second author has been supported by
the Hermann-Minkowski-Minerva Center for Geometry at the Tel Aviv
University and by the German-Israeli
Foundation grant no. 1174-197.6/2011. Main ideas behind this work have appeared during
the visit of the second author to the Higher School of Economics (Moscow),
and a part of this work has been done during
the second author's visits to the Technische Universit\"at Kaiserslautern and
to the Centre Interfacultaire Bernoulli at EPFL, Lausanne. We are grateful to
these institutions for support and excellent working conditions.
We thank Gert-Martin Greuel for useful discussions and Ragni Piene for
important remarks on the preliminary version of the paper
and providing us with references on the subject.}

\section{Equiclassical families of plane curves and duality}\label{sec2}

The discriminant $\Delta_d$ in the linear system $|{\mathcal O}_{\PP^2}(d)|$ of plane curves of
a given degree $d$ admits natural stratifications into families, parameterizing singular
curves with specific properties. In particular, the part
$\Delta_d^{irr}\subset\Delta_d$ formed by reduced, irreducible curves
splits into the disjoint union of Severi varieties $V_{d,g}$ that parameterize curves of degree
$d$ and genus $g$. These are the most important strata of the discriminant, since they are tightly related to
the Gromov-Witten theory
and geometry of moduli spaces of curves. Severi varieties are well-studied: in particular (see \cite{S2,H5}),
for any $d$ and $g$ satisfying
$$0\le g\le\frac{(d-1)(d-2)}{2}\ ,$$
$V_{d,g}$ is a non-empty quasiprojective irreducible variety of dimension
$3d+g-1$; furthermore, a generic member of $V_{d,g}$ is a nodal curve
\cite{AC1,H5,Z2} (a similar statement for other rational surfaces can be found in
\cite{DS}). More thorough stratification consists of
equisingular families of curves that parameterize curves of a given degree
having singular points of prescribed topological or analytic types. Except for small degrees or for curves with only
ordinary nodes, no complete description of such a stratification is known, even the non-emptiness question
on families of curves with ordinary nodes and cusps is not completely answered
(we refer the reader to \cite{GLS1} for a survey and references on this subject).
In the present paper, we address the regularity properties of the
equiclassical stratification of $\Delta^{irr}_d$, which is intermediate between
the Severi and the equisingular stratification. It consists of the
equiclassical families $V_{d,g,c}$ that parameterize curves of
degree $d$, genus $g$, and class $c$ (degree of the dual curve).
A remarkable property of the equiclassical stratification is a series of coupling duality isomorphisms.

\begin{lemma}\label{ldual} The duality transformation of plane curves defines isomorphisms
\begin{equation}
\cD_{d,g,c}:V_{d,g,c}\overset{\simeq}{\longrightarrow}V_{c,g,d},\quad
\cD_{c,g,d}:V_{c,g,d}\overset{\simeq}{\longrightarrow}V_{d,g,c}
\label{eec4}\end{equation} for all $d,c\ge2$, $g\ge0$.
\end{lemma}

{\bf Proof.}
We only
explain that $\cD_{d,g,c}$ and $\cD_{c,g,d}$ are morphisms between quasiprojective varieties. Let $U_d\subset|{\mathcal O}_{\PP^2}(d)|$
be the set of smooth curves of degree $d$. The duality transformation $\cD_d:U_d\to|{\mathcal O}_{\PP^2}(d(d-1))|$
is a well-defined morphism, which can explicitly be written as a homogeneous polynomial map in coefficients of curve equations.
The same polynomial formulas correctly extend to $U_d\cup\Delta^{irr}$ so that, for each curve
$C\in V_{d,g,c}\subset\Delta^{irr}$ with $c<d(d-1)$, the image $\cD_d(C)$ is a reducible curve
containing a component $C'=\cD_{d,g,c}(C)\in V_{c,g,d}$ and $d(d-1)-c$ lines (these
lines are multiple tangents to $C'$, counted with appropriate multiplicities). Hence,
$\cD_d$ restricts to $V_{d,g,c}$ as a well-defined morphism $\cD_{d,g,c}:V_{d,g,c}\to
V_{c,g,d}$. \proofend

Equiclassical families have been studied by Diaz and Harris \cite{DH}.
Answering a question by W. Fulton, they proved

\begin{theorem}[\cite{DH}, Theorem 1.2]\label{tdh}
If $V_{d,g,c}\ne\emptyset$ and
\begin{equation}c\ge2g-1\ ,\label{eec1}\end{equation} then a generic member of any component of
$V_{d,g,c}$ is a curve with nodes and cusps.
\end{theorem}

Later on it was slightly improved by the second author:

\begin{theorem}[\cite{Sh1}, Theorem 1.1]\label{shu}
If $V_{d,g,c}\ne\emptyset$ and
\begin{equation}c\ge2g-d+2\ ,\label{eec2}\end{equation}
then a generic member of any component of
$V_{d,g,c}$ is a curve with nodes and cusps.
\end{theorem}

It is not difficult to show (see \cite{DH,Sh1} or Theorem \ref{tlc1} below) that,
under either of conditions \eqref{eec1}
or \eqref{eec2}, $V_{d,g,c}$ is equidimensional of (expected) dimension
\begin{equation}\dim V_{d,g,c}=\dim_{exp}V_{d,g,c}=d-g+c+1\ .\label{eec5}\end{equation}

Our aim is to study geometry of equiclassical families in more detail.
We will give a series of sufficient numerical conditions for the non-emptiness,
expected dimension, irreducibility, and the
nodal-cuspidal generic member of (any component of) an equiclassical family $V_{d,g,c}$, as well as for a
regular adjacency of equiclassical strata:
\begin{itemize}\item in Section \ref{lc} we mainly revise known results (including those
in \cite{DH,K1,K2,Sh1}) and produce
``linear" conditions, where the number of virtual (in the sense of
 \cite{Ku}) cusps is at most linear in the degree,
\item in Section \ref{sec1}, using methods of \cite{GLS2} and \cite{Sh2}, we derive ``quadratic" conditions,
where the number of (virtual) cusps is bounded from above by a quadratic function of the degree,
\item furthermore, we intertwine both, linear and quadratic conditions
with the duality isomorphisms \eqref{eec4}, which results
in additional regularity conditions, not covered by any other criterion.
\end{itemize}

\subsection{Local structure of equiclassical families}\label{local}

Let $C\in V_{d,g,c}$, $z\in\Sing(C)$. The $\delta$- and $\kap$-invariants of the germ $(C,z)$ are
(see \cite[Definition 3.12]{DH} or \cite[Section I.3.4]{GLS3} for details):
\begin{equation}
\kap(C,z)=(C\cdot C_p)_z,\quad\delta(C,z)=\frac{1}{2}\left(\kap(C,z)-\mt(C,z)+r(C,z)\right)\
, \label{eec12}\end{equation} where $C_p$ denotes a generic polar of $C$, $(*\cdot*)_z$ the
intersection number at the point $z$, $\mt(C,z)=(C\cdot L)_z$ the multiplicity of the singular point
$z$ of $C$, $L$ being a generic straight line passing through $z$, and,
finally, $r(C,z)$ denotes the number of irreducible components of
the germ $(C,z)$.

Furthermore, given
a (local) equation $f(x,y)=0$ of the germ $(C,z)$, we identify the semiuniversal deformation base $B_z$
of the germ $(C,z)$ with ${\mathcal O}_{\PP^2,z}/J(C,z)$, where
$J(C,z)=\langle f,\frac{\partial f}{\partial x},\frac{\partial f}{\partial y}\rangle\subset{\mathcal O}_{\PP^2,z}$
is the Tjurina ideal (cf. \cite[Section 5]{DH}). By \cite{D} and \cite[Definition 4.3 and Theorem 5.5]{DH},
the equiclassical locus
$EC_z\subset B_z$ (which parameterizes $\delta$- and $\kap$-constant deformations)
is an irreducible analytic variety germ, which has codimension $\kap(C,z)-\delta(C.z)$ with the tangent cone
supported at the linear space $J^{ec}(C,z)/J(C,z)$, where $J^{ec}(C,z)\subset{\mathcal O}_{C,z}$
is the ideal defined by
$$J^{ec}(C,z)=\left\{\varphi\in{\mathcal O}_{\PP^2,z}
\ \Bigg|\ \begin{array}{l}\ord\varphi\big|_P\ge(C_p\cdot P)\;\ \text{for each}\\
\text{irreducible component}\ P\subset(C,z)\end{array}
\right\}\ .$$
We have a natural map
$$\Phi_C:\left(|{\mathcal O}_{\PP^2}(d)|,C\right)\to\prod_{z\in\Sing(C)} B_z\ ,$$ and the germ of $V_{d,g,c}$ at $C$
can be viewed as
$$\left(V_{d,g,c},C\right)=\Phi^{-1}\left(\prod_{z\in\Sing(C)}EC_z\right)\ .$$
Notice that $EC_z$ and hence $V_{d,g,c}$ may not be smooth (see \cite[Theorem 28 and Lemma 29]{D}).
We, however, can speak about the following property:

\begin{definition}\label{dec2}
We say that the germ of the family $V_{d,g,c}$ at
$C\in V_{d,g,c}$ is a {\bf locally transverse intersection} (shortly {\bf LTI}), if
$\Ima\Phi_C$ intersects transversally with $\prod_{z\in\Sing(C)}EC_z$ in $\prod_{z\in\Sing(C)}B_z$.
Respectively, $V_{d,g,c}$ is called LTI, if all germs $(V_{d,g,c},C)$ are LTI for all $C\in V_{d,g,c}$.
\end{definition}

It is evident that the duality isomorphisms \eqref{eec4} preserve the expected dimension of
equiclassical families (cf. formula \eqref{eec5}). We can say even more:

\begin{theorem}\label{tec3}
(1) Let $V_{d,g,c}$ be LTI. Then $V_{d,g,c}$ has expected dimension $c-g+d+1$, and satisfies the following
incidence relations:
\begin{equation}\begin{cases}V_{d,g,c}\subset\overline V_{d,g,c+1},\quad & \text{if}\quad 2d-2-c+2g\\
&\quad=\sum_{z\in\Sing(C)}(\kappa(C,z)-2\delta(C,z))>0,\\
V_{d,g,c}\subset\overline V_{d,g+1,c+2},\quad & \text{if}\quad c-3g+\frac{1}{2}(d^2-7d+6)\\
&\quad =\sum_{z\in\Sing(C,z)}(3\delta(C,z)-\kappa(C,z))>0,
\end{cases}\label{eec30}\end{equation} $C$ being any curve in $V_{d,g,c}$.

(2) The family $V_{d,g,c}$ is LTI if and only if the family $V_{c,g,d}$ is LTI.
\end{theorem}

\begin{remark}
Note that the first inclusion in \eqref{eec30} means deforming a cusp to a node, while the second inclusion means
smoothing out a node of a generic curve in $V_{d,g,c}$.
\end{remark}

{\bf Proof.} (1) For the first statement, we notice that the LTI property yields
$$\dim V_{d,g,c}=\frac{d(d+3)}{2}-\sum_{z\in\Sing(C)}(\kap(C,z)-\delta(C,z))
=d-g+c+1\ ,$$ where the latter equality follows from the Pl\"ucker and genus formulas
\begin{equation}c=d(d-1)-\kap(C),\quad g=\frac{(d-1)(d-2)}{2}-\delta(C)\ ,
\label{eec21}\end{equation} where $$\kap(C)=\sum_{z\in\Sing(C)}\kap(C,z),\quad \delta(C)=\sum_{z\in
\Sing(C)}\delta(C,z)\ .$$
Now pick $C\in V_{d,g,c}$ and show that $C\in\overline V_{d,g,c+1}$
or $C\in\overline V_{d,g+1,c+2}$ according as $2d-2-c+2g>0$ or
$c-3g-\frac{1}{2}(d^2-7d+6)>0$. Since
$$2d-2-c+2g=\kap-2\delta,\quad c-3g+\frac{d^2-7d+6}{2}=3\delta-\kap\ ,$$
we can interpret these numbers
as follows. By \cite[Theorem 1.1]{DH}, a generic germ of
an equiclassical locus $EC_z$, where $z\in\Sing(C)$, $C\in V_{d,g,c}$, has $\kap(C,z)-2\delta(C,z)$
cusps and $3\delta(C,z)-\kap(C,z)$ nodes as its only singularities. Hence, $\kap-2\delta$
and $3\delta-\kap$ are
the total numbers of cusps, resp. nodes over generic members of $EC_z$, $z\in\Sing(C)$.
If $\kap-2\delta>0$, then there is $z\in\Sing(C)$ such that $\kap(C,z)-2\delta(C,z)>0$,
and hence $EC_z\subset\overline{EC}'_z$, where $EC'_z\subset B_z$ is the locus
formed by the germs having $\kap(C,z)-2\delta(C,z)-1$ cusps and $3\delta(C,z)-
\kap(C,z)+1$ nodes\footnote{Over a versal deformation base, all singularities can be
independently deformed or preserved.}. It follows from the LTI property that
$C\in\overline V_{d,g,c+1}$. If $3\delta-\kap>0$, then there is
$z\in\Sing(C)$ such that $3\delta(C,z)-\kap(C,z)>0$,
and hence $EC_z\subset\overline{EC}''_z$, where $EC''_z\subset B_z$ is the locus formed by the
germs having $\kap(C,z)-2\delta(C,z)$ cusps and $3\delta(C,z)-\kap(C,z)-1$ nodes.
Again the LTI property yields that $C\in\overline V_{d,g+1,c+2}$.

(2) For the second statement, for each point $z\in\Sing(C)$, introduce the
zero-dimensional scheme $Z^{ec}_{C,z}\subset\PP^2$ supported at $z$ and defined by the ideal
$J^{ec}(C,z)$. Set $Z^{ec}(C)=\bigcup_{z\in\Sing(C)}Z^{ec}_{C,z}$. By
\cite[Lemma 5.3]{DH}, $\deg Z^{ec}(C)=\kap-\delta$. The LTI property
can equivalently be expressed as the $h^1$-vanishing
\begin{equation}H^1(\PP^2,{\mathcal J}_{Z^{ec}(C)}(d))=0\ ,\label{eec13}\end{equation}
where ${\mathcal J}_{Z^{ec}(C)}$ is the ideal sheaf of the scheme $Z^{ec}(C)$.
We rewrite condition \eqref{eec13}
as follows. Let $F(x_0,x_1,x_2)=0$ be a (homogeneous) equation of a curve $C\in V_{d,g,c}$.
Choose a generic point $p=(p_0,p_1,p_2)$ and a generic line $L_q=q_0x_0+q_1x_1+q_2x_2=0$ in the plane so that
\begin{itemize}\item the polar curve $C_p=\{p_0F_{x_0}+p_1F_{x_1}+p_2F_{x_2}=0\}$ intersects $C$ at each
points $z_i\in\Sing(C)$ with multiplicity $(C_p\cdot C)_{z_i}=\kap(C,z_i)$, $i=1,...,r$,
and in a set $S_p\subset C\setminus\Sing(C)$ with the total multiplicity $c$;
\item the set
$T_q=L_q\cap C$ is disjoint from $\Sing(C)\cup S_p$, and $(L_q\cdot C)_{T_q}=d$. \end{itemize}
Then (cf. \cite[Proof of Lemma 5.14]{DH}, \cite[Proof of Lemma 2.1]{Sh1})
\begin{equation}
H^1(\PP^2,{\mathcal J}_{Z^{ec}(C)}(d))=H^1(\hat C,{\mathcal O}_{\widehat C}(\nu^*(S_p)+\nu^*(T_q)))\ ,
\label{eec23}\end{equation} where $\nu:\widehat C\to C$ is the normalization. Let $C^\vee\in V_{c,g,d}$ be the dual curve,
\mbox{$\nu^\vee:\widehat C\to C^\vee$} the dual normalization, and let $C^\vee$ be given
by a homogeneous equation $F^\vee(y_0,y_1,y_2)=0$. Then the polar $C^\vee_q=
\{q_0F^\vee_{y_0}+q_1F^\vee_{y_1}+q_2F^\vee_{y_2}=0\}$ intersects $C^\vee$ in $\Sing(C^\vee)$ and in a set
$S^\vee_q$ dual to $T_q\subset C$, whereas the line $L^\vee_p=\{p_0y_0+p_1y_1+p_2y_2=0\}$ intersects
$C^\vee$ in a set $T^\vee_p\subset C^\vee\setminus(\Sing(C^\vee)\cup S^\vee_q)$ dual to $S_p\subset C$.
Hence $\nu^*(S_p)=(\nu^\vee)^*(T^\vee_p)$, $\nu^*(T_q)=
(\nu^\vee)^*(S^\vee_q)$, which completes the proof in view of
$$H^1(\PP^2,{\mathcal J}_{Z^{ec}(C^\vee)}(c))=H^1(\hat C,{\mathcal O}_{\widehat C}
((\nu^\vee)^*(S^\vee_q)+(\nu^\vee)^*(T^\vee_p)))\ .\quad\quad\quad\quad\text{\proofend}$$

\begin{remark}\label{rec4}
One can see from the proof of Theorem \ref{tec3} that the values
$2d-2-c+2g$ and $c-3g+\frac{1}{2}(d^2-7d+6)$ are always nonnegative, and their sum
$\frac{1}{2}(d^2-3d+2)-g$ is positive as long as $V_{d,g,c}$ contains singular curves.
\end{remark}

\begin{definition}\label{dec1}
We call a family $V_{d,g,c}$ {\bf locally regular} if it is LTI,
has expected dimension and satisfies
the incidence relations \eqref{eec30}.
\end{definition}

A simple well-known example of locally regular families are Severi varieties
$V_{d,g}$ (here $c=2d-2+2g$), see \cite{S2}.

\subsection{Linear regularity conditions}\label{lc} In the following theorem we
present some elementary or known regularity conditions.

\begin{theorem}\label{tlc1}
(1) The family $V_{d,g,c}$, $d,c\ge2$, is non-empty, if
\begin{equation}c\ge2g+\left[\frac{d+1}{2}\right]+1
\quad\left(\text{equivalently}\quad
\kappa-2\delta\le d+\left[\frac{d}{2}\right]-3\right)\ .\label{eec26}\end{equation}

(2) Let $V_{d,g,c}\ne\emptyset$, $d,c\ge2$.
\begin{enumerate}\item[(2i)] If
\begin{equation}c-2g+d\ge-1\quad(\text{equivalently}\quad
\kappa-2\delta\le3d-1)\ ,\label{eec22}\end{equation}
then $V_{d,g,c}$ is locally regular.
\item[(2ii)] If
\begin{equation}c-2g+d\ge2\quad(\text{equivalently}
\quad\kappa-2\delta\le3d-4)\ ,\label{eec29}\end{equation} then a generic member of any
component of $V_{d,g,c}$ is a curve with nodes and cusps.
\item[(2iii)] If
\begin{equation}\begin{cases}\text{either}\ &c\ge2g+2d-5\quad(\text{equivalently}\
\kappa-2\delta\le3),\\
\text{or}\ &c\ge3g+\frac{3d-5}{2}\quad\ \ (\text{equivalently}\
3\delta-\kappa\ge\frac{d^2-4d+1}{2}),\\
\quad\text{or}\ &c\ge g+\frac{d^2-2d-1}{2}\quad(\text{equivalently}
\ \kappa-\delta\le\frac{3d-1}{2}),\end{cases}
\label{eec24}\end{equation}
then $V_{d,g,c}$ is irreducible.
\end{enumerate}
\end{theorem}

{\bf Proof.} Under condition \eqref{eec26}, we will construct curves with $n$ nodes and $k$ cusps
realizing the required values of $d,g,c$. We can rewrite \eqref{eec26} as
$k\le\left[\frac{3(d-2)}{2}\right]$. It is classically known (see \cite{S3,S4,S2}, and
also \cite[Corollary 6.3]{GK} for a modern treatment), that
if $k<3d$, then one can smooth out prescribed nodes and cusps of such a nodal-cuspidal curve and deform
prescribed cusps into nodes. Hence, it is sufficient to construct the following
nodal-cuspidal curves: (a) a rational curve of degree $d=2s$ with $k=3s-3$, (b) a rational curve of degree
$d=2s+1$ with $k=3s-2$. By Pl\'ucker formulas, the first curve is dual to a nodal rational curve
of degree $c=s+1$, whereas the second curve is dual to a rational curve of
degree $c=s+2$ having one cusp and $(s^2+s-2)/2$ nodes.
The existence of the latter curve is left to the reader as an elementary exercise.

Claim (2ii) is a copy of Theorem \ref{shu}, included here for completeness and comparison with
other criteria. Claim (2i) (which basically coincides with
\cite[Corollary 5.2]{DH}) holds, since \eqref{eec22} implies
\eqref{eec13}: indeed, in view of \eqref{eec23} it is enough to show that
$$2g-2<\deg(\nu^*(S_p)+\nu^*(T_q))=c+d\ ,$$ which is equivalent to \eqref{eec22}.
In Claim (2iii) we can suppose that $d,c\ge3$, and then each inequality in \eqref{eec24}
yields \eqref{eec2}: it is evident for the first two inequalities, and the third one
combined with $g\le2g-d+1$ negative to \eqref{eec2} results in the impossible
relation $g\ge\frac{1}{2}(d^2-1)$. Hence, by Theorem \ref{shu},
we are left with the irreducibility problem for families
$V_d(nA_1,kA_2)$ of irreducible curves of degree $d$ with $n$ nodes and $k$ cusps. Each of the following
conditions is sufficient for the irreducibility of $V_d(nA_1,kA_2)$:
\begin{itemize}\item $k\le3$ (see \cite{K1}, the missing cases of $k=3$ and $d=5$ or $6$ are
covered in \cite[Theorem 2.4]{AD}, \cite{Deg1}, \cite[Section II.7.3]{Deg});
\item $k\le\frac{1}{2}(d+1)-g$ (see \cite{K2});
\item $2n+4k<3d$ (see \cite[Corollary IV.6.2]{GLS4}).
\end{itemize}
In terms of $d,g$, and $c$ they specialize to the form \eqref{eec24}.
\proofend

\begin{remark}\label{rec3}
The bound \eqref{eec26} is optimal for rational curves: a larger number of cusps
simply is not possible by Pl\"ucker formulas as one can see in the proof of Theorem \ref{tlc1}.
\end{remark}

Now we combine Theorem \ref{tlc1} with isomorphisms \eqref{eec4}. Notice, first, that
linear conditions remain linear as long as $k\le3d$: Indeed, it follows from Pl\"ucker formulas that
$3d-k=3c-k^\vee$, where $k$, resp. $k^\vee$, is the number of (virtual) cusps
of the original, resp. dual curve. Second, the dualization of Claims (1), (2i), and (2ii)
in Theorem \ref{tlc1}
does not give anything new, whereas the dual to Claim (2iii) is a {\it new} statement:

\begin{theorem}\label{tlc2}
(1) If $V_{d,g,c}\ne\emptyset$ and
\begin{equation}\text{either}\ d\ge2g+2c-5,\quad\text{or}\ d\ge3g+\frac{3c-5}{2}\ ,
\quad\text{or}\ d\ge g+\frac{c^2-2c-1}{2}\ ,
\label{eec27}\end{equation}
then $V_{d,g,c}$ is irreducible.

(2) A non-empty family $V_d(nA_1,kA_2)$ is irreducible if
\begin{equation}\begin{cases}\text{either}\quad & 3n+4k\ge\frac{1}{2}(3d^2-6d-3)\ ,\\
\text{or}\quad &
4n+5k\ge\frac{1}{3}(6d^2-14d+1)\ ,\\
\text{or}\quad & n+k\ge\frac{1}{2}(d^2-5d+(d^2-d-2n-3k-1)^2)\ .
\end{cases}\label{eec28}\end{equation}
\end{theorem}

{\bf Proof.} Conditions \eqref{eec27} are obtained from \eqref{eec24} by exchange of $d$ and $c$.
Furthermore, rewriting \eqref{eec27} in the form
$$c-2g+d\ge3c-5,\quad c-2g+d\ge g+\frac{5c-5}{2},\quad
c-2g+d\ge\left(\frac{c^2-3c+2}{2}-g\right)+\frac{2c-3}{2},$$ respectively, we see that each of the inequalities
\eqref{eec27} yields \eqref{eec29}, and hence we get irreducibility criteria for families of curves with
nodes and cusps in Claim (2): via Pl\"ucker formulas, inequalities \eqref{eec27}
are converted into \eqref{eec28}.\proofend

\begin{example}\label{exec1}
To show that Theorem \ref{tlc2} indeed contains a new information as compared with Theorem \ref{tlc1},
we consider the case of curves of degree $d=7$ (for $d\le6$
the equisingular stratification of the discriminant is completely described
\cite{AD,Deg1,Deg}). So, Theorem \ref{tlc1}(2iii) implies that $V_7(nA_1,kA_2)$ is irreducible if
$$\text{either}\ k\le3,\quad\text{or}\ k=4,\ n=11,
\quad\text{or}\ k=4,\ n\le3
\ .$$ In addition to these cases, Theorem \ref{tlc2} implies the irreducibility of
$V_7(nA_1,kA_2)$ for
$$(k,n)=(6,9),\ (7,8),\ \text{and}\ (10,4)\ .$$
All these families are non-empty. One can easily find similar examples in higher degrees.
\end{example}

\subsection{Quadratic regularity conditions}\label{sec1}

To formulate conditions, where the number of (virtual) cusps may
depend quadratically on the degree, we replace the parameters $d,g,c$ by
$d,\delta,\kap$ and use the notation
$V_d^{\delta,\kap}=V_{d,g,c}$, where $\kap=d(d-1)-c$, $\delta=(d-1)(d-2)/2-g$ (that is
$\delta=\delta(C)$ and $\kap=\kap(C)$ for any
curve $C\in V_d^{\delta,\kap}=V_{d,g,c}$). For example, in terms of $\kap$ and $\delta$, condition
\eqref{eec1}, resp. \eqref{eec2} reads
\begin{equation}\kap-2\delta\le2d-1,\quad\text{resp.}\quad
\kap-2\delta\le3d-4\ .\label{eec3}\end{equation}

Due to the irregular behavior of families of plane curves with nodes and cusps
(for instance, it is not known how many cusps a curve of degree $d\ge9$ may have, cf. \cite{CPS}),
it is not possible to completely describe non-empty equiclassical
families. Here we present the following partial result:

\begin{theorem}\label{tec2}
Let $d\ge3$, $\kap,\delta>0$.

(1) If $V_d^{\delta,\kap}\ne\emptyset$, then
\begin{equation}2\delta\le\kap\le3\delta,\quad\delta\le\frac{(d-1)(d-2)}{2}\
.\label{eec10}\end{equation}

(2) If $d\le4$, then \eqref{eec10} is sufficient for the non-emptiness of $V_d^{\delta,\kap}$.

If $d\ge5$, then the conditions \eqref{eec10} and
\begin{equation}\kap-\delta\le\frac{d^2-4d+6}{2}\label{eec11}\end{equation}
are sufficient for the non-emptiness of $V_d^{\delta,\kap}$.
\end{theorem}

{\bf Proof.} (1) The first inequality in \eqref{eec10} follows from the second relation in \eqref{eec12} and
from the relation
$\delta(C,z)\ge\frac{1}{2}\mt(C,z)(\mt(C,z)-1)$, which follows from \eqref{eec12} too.

(2) By \cite[Corollary 1.3]{Sh1}, for $d\le10$ a generic member of any component of
a non-empty family $V_d^{\delta,\kap}$ is a curve with nodes and cusps. Thus, the case of
$d\le4$ reduces to the study of curves with nodes and cusps, which are classically known.
For $d\ge5$ we apply \cite[Theorem 4.1]{Sh2} (suitably specified in
\cite[Theorem IV.5.4(ii)]{GLS4}), which states that the inequality
$n+2k\le\frac{1}{2}(d^2-4d+6)$ is sufficient for the non-emptiness of
$V_d(nA_1,kA_2)$. In terms of $\delta=n+k$ and $\kap=2n+3k$, the latter inequality
reads \eqref{eec11}.
\proofend

\begin{remark}\label{rec1}
The non-emptiness condition \eqref{eec11} is asymptotically sharp for equiclassical families
having expected dimension, since their codimension satisfies
$$\kap-\delta\le\dim|{\mathcal O}_{\PP^2}(d)|=\frac{d(d+3)}{2}\ .$$
\end{remark}

\begin{theorem}\label{tqc1}
Let $d\ge3$, $\kap,\del\ge0$, and let $V_d^{\delta,\kap}\ne\emptyset$.

(1) If
\begin{equation}5\kap-6\delta\le(d+3)^2\ ,\label{eec7}\end{equation}
then $V_d^{\delta,\kap}$ is locally regular and
$$\dim V_d^{\delta,\kap}=\frac{d(d+3)}{2}-\kap+\delta\ .$$

(2) If
\begin{equation}5\kap-6\delta\le d^2+6d-3\ ,\label{eec8}\end{equation}
then a generic member of any component of $V_d^{\delta,\kap}$ is a curve
with $3\delta-\kap$ nodes and $\kap-2\delta$ cusps as its only singularities.

(3) If
\begin{equation}\frac{11}{2}\;\kap+\frac{3}{2}\;\delta<d^2\ ,\label{eec9}\end{equation}
then $V_d^{\delta,\kap}$ is irreducible.
\end{theorem}

\begin{remark}\label{rec2}
If a generic member of an irreducible component of $V_d^{\delta,\kap}$ is a curve with
$n$ nodes and $k$ cusps, then inequality \eqref{eec7} reads
$4n+9k\le(d+3)^3$, the best known quadratic sufficient condition for
the smoothness and expected dimension of the family of plane curves of degree $d$ with $n$ nodes and $k$ cusps
\cite[Corollary 4.4]{GLS2}.
\end{remark}

{\bf Proof of Theorem \ref{tqc1}.} (1) The LTI property is equivalent to the $h^1$-vanishing
\eqref{eec13}.
Since $J^{ec}(C,z_i)\supset J(C,z_i)$ for all $i=1,...,r$, we can apply the $h^1$-vanishing criterion
from \cite[Proposition 4.1]{GLS2} sufficient for \eqref{eec13}:
\begin{equation}\sum_{i=1}^r\gamma(C;Z^{ec}_{C,z_i})\le(d+3)^3\ ,
\label{eec16}\end{equation} where the $\gamma$-invariant can be computed as follows: for an
irreducible zero-dimensional scheme
$Z\subset C$ supported at a point $z\in\Sing(C)$ and such that its defining ideal satisfies $I_Z\supset J(C,z)$,
\begin{equation}\gamma(C;Z)=\max\left\{4\deg Z,\ \max_D\frac{(D\cdot C)_z^2}{(D\cdot C)_z-\deg Z}\right\}\ ,
\label{eec14}\end{equation}
where $D$ runs over all plane curve germs defined by elements $g\in I_Z$ such that $(D\cdot C)_z\le2\deg Z$
(see \cite[Formula (4.1)]{GLS2} and \cite[Definition I.2.25]{GLS4}). Let us show that
\begin{equation}\gamma(C;Z^{ec}_{C,z_i})=\frac{\kap(C,z_i)^2}{\delta(C,z_i)}\ .
\label{eec15}\end{equation} Indeed, by definition of
$Z^{ec}_{C,z_i}$, always $(D\cdot C)_z\ge\kap(C,z_i)$, and hence the maximal value of the right-hand side of
\eqref{eec14} is achieved for $(D\cdot C)_z=\kap(C,z_i)$, which together with
$\deg Z^{ec}_{C,z_i}=\kap(C,z_i)-\delta(C,z_i)$ yields \eqref{eec15}. Observing
that $a^2/b\le5a-6b$ as long as $0<2b\le a\le 3b$, we derive the sufficient condition
\eqref{eec7} for \eqref{eec16}.

\smallskip

(2) Let $V$ be an irreducible component of $V_d^{\delta,\kap}$, $C\in V$ its generic member
having singular points $z_1,...,z_r$ of topological types $S_1,...,S_r$, respectively. Then
the germ of $V$ at $C$ is the germ at $C$ of an equisingular family $V_d(S_1,...,S_r)$.
Suppose that $z_1$ is not a node or a cusp.
By \cite[Pages 144, 158, 164]{W} (see also \cite[Theorem 2.25]{DH}), the equisingular
stratum $ES_1\subset B_1$ is smooth and its tangent space is $J^{es}(C,z_1)/J(C,z_1)$,
where $J^{es}(C,z_1)\subset{\mathcal O}_{\PP^2,z_1}$ is the equisingular ideal. Furthermore,
by \cite[Proposition 5.10 and Theorem 5.12]{DH}, $J^{ec}(C,z_1)\subsetneq J^{es}(C,z_1)$.
Then there exists an intermediate ideal $J^{ec}(C,z_1)\subsetneq J\subset J^{es}(C,z_1)$ such that
$\dim J/J^{ec}(C,z_1)=1$. Introduce the zero-dimensional scheme $Z'=Z_1\cup\bigcup_{i\ge2}Z^{ec}_{C,z_i}$
with $Z_1$ defined at $z_1$ by the ideal $J$. Then
\begin{equation}\dim_CV\le h^0(\PP^2,{\mathcal J}_{Z'}(d))-1\ .\label{eec17}\end{equation}
Using \eqref{eec8}, we will prove that \begin{equation}H^1(\PP^2,
{\mathcal J}_{Z'}(d))=0\ ,\label{eec18}\end{equation} which in turn
will imply
\begin{eqnarray}\dim_C V&\le& h^0(\PP^2,{\mathcal J}_{Z'}(d))-1=\frac{d(d+3)}{2}-\deg Z'\nonumber\\
&=&\frac{d(d+3)}{2}-\deg Z^{ec}(C)-1=\frac{d(d+3)}{2}-\kap+\delta-1\ ,
\nonumber\end{eqnarray} a contradiction to
Claim (1). To establish \eqref{eec18}, we again apply \cite[Proposition 4.1]{GLS2} and verify
the inequality
\begin{equation}\gamma(C;Z_1)+\sum_{i=2}^r\gamma(C;Z^{ec}_{C,z_i})\le(d+3)^2\ .\label{eec19}\end{equation}
Similarly to the above computation of $\gamma(Z^{ec}_{C,z_i})$, we obtain
$$\gamma(C;Z_1)\le\frac{\kap(C,z_1)^2}{\kap(C,z_1)-\deg Z_1}=
\frac{\kap(C,z_1)^2}{\delta(C,z_1)-1}\ .$$ It is an easy exercise to show (using formulas
\eqref{eec12}) that $\delta(C,z_1)\ge2$ and \mbox{$\kap(C,z_1)\le\frac{5}{2}\delta(C,z_1)$}, which yields
$$\frac{\kap(C,z_1)^2}{\delta(C,z_1)-1}\le5\kap(C,z_1)-6\delta(C,z_1)+12\ .$$ Joining this with
the inequality $\kap(C,z_i)^2/\delta(C,z_i)\le5\kap(C,z_i)-6\delta(C,z_i)$ obtained in Step (1),
we see that \eqref{eec8} implies \eqref{eec19}, and hence \eqref{eec18}, completing the proof of Claim (2).

\smallskip

(3) Observe that \eqref{eec9} yields \eqref{eec8}. Hence the irreducibility of
$V_d^{\delta,\kap}$ will follow from the
irreducibility of the family of plane curves of degree $d$ with $n=3\delta-\kap$ nodes and
$k=\kap-2\delta$ cusps. By \cite[Theorem 2]{GLS2}, $V_d(nA_1,kA_2)$ is irreducible when
$25n/2+18k<d^2$, which is equivalent to \eqref{eec9}.
\proofend

To apply duality isomorphisms \eqref{eec4}, we go back to the parameters $d,g,c$:

\begin{theorem}\label{tqc2}
(1) The family $V_{d,g,c}$ is non-empty if $g\le\frac{1}{2}(d-1)(d-2)$ and
$$\text{either}\ c\ge g+\frac{5d-8}{2}\;,\quad\text{or}\ d\ge g+\frac{5c-8}{2}\ .$$

(2) Let $V_{d,g,c}\ne\emptyset$. Then
\begin{enumerate}\item[(2i)] $V_{d,g,c}$ is locally regular if
$$\text{either}\ 5c\ge6g+d^2-2d-15,\quad\text{or}\ 5d\ge6g+c^2-2c-15\ ;$$
\item[(2ii)] $V_{d,g,c}$ is irreducible if
$$\text{either}\ 11c+3g>\frac{21d^2-31d+6}{2}\;,\quad
\text{or}\ 11d+3g>\frac{21c^2-31c+6}{2}\ .$$
\end{enumerate}
\end{theorem}

Here each of the ``either" inequality is a translation of the corresponding condition in Theorem \ref{tqc1},
and each of the ``or" inequality is the dual. We notice that the nodal-cuspidal property is
a priori not preserved by the duality.

\begin{example}\label{exec2}
We again illustrate the novelty of the dual criteria. Consider the family
$V_d(kA_2)$ with large $d$ and $k=[(d+3)^2/9]$. It is non-empty by
Theorem \ref{tec2}(2) and is LTI of expected dimension by Theorem
\ref{tqc1}(1) (cf. Remark \ref{rec2}). The dual curves have degree and genus
$$d^\vee=c=d(d-1)-3k\sim\frac{2}{3}\;d^2,\quad g=\frac{(d-1)(d-2)}{2}-k\sim\frac{7}{18}\;d^2\ ,$$ which correspond to
the following numbers of virtual nodes and cusps
$$n^\vee\sim\frac{2}{9}\;d^4,\quad k^\vee\sim\frac{19}{9}\;d^2\ ,$$
but these values do not satisfy the other known LTI/expected dimension conditions: neither
$k^\vee<3d^\vee$, nor $4n^\vee+9k^\vee<(d^\vee+3)^2$.
\end{example}

\section{Duality of spacial curves and associated maps}\label{sec3}

Duality for non-degenerate curves in $\PP^n$, $n\ge3$, was discovered by R. Piene
\cite[\S 5]{P}. It differs from the general duality between subvarieties of $\PP^n$, and it always takes
curve to curve. It, furthermore, can be extended to curves in Grassmannians \cite{Per}. Relations
between projective curves and curves in Grassmannians have been treated in \cite{GH1}. We present here
some basic material from \cite{GH,GH1,Per,P} on spacial and associated curves, on their duality, and on associated maps between families of
curves in projective spaces and Grassmannians. For the reader's convenience we provide elementary proofs,
which can be traced back to H. Weyl \cite{We} and which require only an undergraduate calculus. Finally, we formulate few problems on the geometry of
equiclassical families of spacial curves and maps between them.

\smallskip

{\bf Associated and dual curves.}
We consider curves in the (complex) projective spaces $\PP^n$, $n\ge2$,
as maps of Riemann surfaces to $\PP^n$. Let ${\mathcal M}_{g,0}(\PP^n,d)$, $d\ge1$, be the moduli space of stable
maps $f:S\to\PP^n$ of Riemann surfaces $S$ of genus $g\ge0$ such that
$[\nu_*(S)]=d[\PP^1]\in H_2(\PP^n)$. It is a quasiprojective variety
(see, for instance, \cite[Theorem 1]{FP}). Denote by
${\mathcal M}_{g,0}^{\nd}(\PP^n,d)\subset{\mathcal M}_{g,0}(\PP^n,d)$ the subspace
of {\bf non-degenerate} curves, that is (isomorphism classes of)
maps $f:S\to\PP^n$ birational onto the image $f(S)$, which is not contained in any
hyperplane in $\PP^n$.

Let a curve $f:S\to\PP^n$, $n\ge2$, be non-degenerate.
Let $s\in S$, and $t$ a regular parameter in the
analytic germ $(S,s)$.
In suitable projective coordinates $X=(x_0,...,x_n)$, the
map $f:(S,s)\to\PP^n$ is given by a vector-function $X(t)=(x_0(t),\dots,x_n(t))$
(cf. \cite[Page 266]{GH}), where
$$x_0(t)=1+O(t),\ x_1(t)=t^{1+\bet_0}(1+O(t)),\ ...,\ x_k(t)=t^{k+\bet_0+...+\bet_{k-1}}(1+O(t)),$$
\begin{equation}...\ ,\ x_n(t)=t^{n+\bet_0+...+\bet_{n-1}}(1+O(t))\ .\label{esc4}\end{equation}
Given $k=0,...,n-1$, the number $\bet_k=\bet_k(f,s)\ge0$ in the parametrization
\eqref{esc4} does not depend on the choice of a local parameter
and is called the {\bf $k$-th ramification index} of $f$ at the point $s\in S$;
the number $\bet_k(C)=\sum_{s\in S}\bet_k(f,s)$ is finite and is called
the {\bf total $k$-th ramification index} of the curve $C$.

A non-degenerate curve $C=[f:S\to\PP^n]$ induces a well-defined
associated curve $C_k=[f_k:S\to Gr(k+1,n+1)]$, which parameterizes osculating $k$-planes to \mbox{$C=C_0$},
for all $k=1,...,n-1$
(see details, for instance, in \cite[Chapter 2, Section 4]{GH}).
In the chosen coordinates \eqref{esc4}, $f_k$ is locally given by the
multivector-function $X(t)\wedge X'(t)\wedge\dots\wedge X^{(k)}(t)$.
Denote $C_k=\As_k(C)$, $k=0,1,...,n-1$, where $C_0=C$.
Observe that $C_{n-1}=\As_{n-1}(C)$ is a curve in the dual projective space
\mbox{$(\PP^n)^*=Gr(n,\C^{n+1})$}. Denote by $(\C^{n+1})^*$ the space dual to $\C^{n+1}$. We have
the following canonical isomorphisms for all $k=1,...,n$:
$$\begin{cases}&*_k:Gr(k,\C^{n+1})\to Gr(n+1-k,(\C^{n+1})^*),\\
&*_k:Gr(k,(\C^{n+1})^*)\to Gr(n+1-k,((\C^{n+1})^*)^*)=Gr(n+1-k,\C^{n+1}),
\end{cases}$$
which send a linear subspace to the space of linear functionals vanishing on it. Clearly,
$*_k*_{n+1-k}=\Id$ for each $k=1,...,n$.

\begin{definition} \label{maindef} Let $C'=[f':S\to\PP^n]$ and $C''=[f'':
S\to(\PP^n)^*]$ be two non-degenerate curves. We say that
$C''$ is $k$-{\bf dual} to $C'$ (resp., $C'$ is $(n-k)$-{\bf dual} to
$C''$) if $*_k(C'_{k-1})=C''_{n-1-k}$, and we say that $C'$ and
$C''$ are {\bf dual} if $C''$ is $k$-dual to $C'$ for each $k=1,\dots,n-1$.\end{definition}

Obviously, if $C^{\prime\prime}$ is $k$-dual to $C'$ for each $k=1,\dots,n-1$
then $C'$ is also $k$-dual to $C^{\prime\prime}$ for each $k=1,\dots, n-1$.
The key fact is the following duality statement:

\begin{theorem}[\cite{P}, Theorem 5.1]\label{tsc}
Let $C$ be a non-degenerate curve in $\PP^n$, $n\ge2$. Then
$C^*=\As_{n-1}(C)$ is non-degenerate
as well, and, moreover, $C$ and $C^*$ are dual curves.
Furthermore,
$\bet_k(C^*)=\bet_{n-1-k}(C)$ for all $k=0,...,n-1$.
\end{theorem}

{\bf Proof.}
Let  $C$ be locally given by a vector-function $X(t)$ in the form \eqref{esc4}.
Recall that (see \cite[Pages 263-264]{GH}), for a non-degenerate curve $C$, at all but finitely many points
$z\in S$, the vectors $X,X',...,X^{(n-1)}$ (computed in the parametrization \eqref{esc4}) are
linearly independent, and $\bet_k(f,s)=0$ for all $k=0,...,n-1$.

Denote $f^*=f_{n-1}$.
The map $f^*:(S,s)\to(\PP^n)^*$ is given by the $n\times n$-minors $y_0(t),...,y_n(t)$ of the
matrix formed by the rows
$X,X',...,X^{(n-1)}$. It is an easy exercise to show  (after reducing the minimal common power
of $t$) that
$$y_0=t^{n+\bet_{n-1}+...+\bet_0}(\lam_0+O(t)),\ y_1=t^{n-1+\bet_{n-1}+...+\bet_1}(\lam_1+O(t)),\ ...\ ,$$
$$y_k=t^{n-k+\bet_{n-1}+...+\bet_k}(\lam_k+O(t)),\ ...\ ,\ y_n=\lam_n+O(t)\ ,$$
where $\lam_0\lam_1...\lam_n\ne0$. It follows that $C^*$ is non-degenerate,
$$\bet_k(f^*,s)=\bet_{n-1-k}(f,s),\quad k=0,...,n-1,\ s\in S\ ,$$
and for a regular
point $s\in S$, the map $f^*:(S,s)\to(\PP^n)^*$ can also be recovered as the unique
(up to proportionality) solution $Y=(y_0,...,y_n)$ to the system
\begin{equation}XY=0=X'Y=...=X^{(k)}Y=...=X^{(n-1)}Y=0\ ,
\label{esc1}\end{equation} where the product of vectors is understood as the sum of products
of the corresponding coordinates. Differentiating equations \eqref{esc1} and combining them to each other,
we obtain that
\begin{equation}X^{(k)}Y^{(l)}=0\quad\text{for all}\quad k,l\ge0,\ k+l\le n-1\ .
\label{esc2}\end{equation} In particular,
\begin{equation} \label{eq} XY=XY'=...=XY^{(k)}=...=XY^{(n-1)}=0\ ,\end{equation} which means that $C=\As_{n-1}(C^*)$.

Furthermore, if $n\ge3$, it follows from \eqref{esc2} that, for each $k=1,...,n-2$, we have
\begin{equation}X^{(i)}Y^{(j)}=0,\quad \text{for all}\quad 0\le i\le k,\ 0\le j\le n-1-k
\ .\label{esc3}\end{equation}
The vectors $Y,Y',...,Y^{(n-1)}$ are linearly independent too, since $\As_{n-1}(C)$ is a non-degenerate curve.
The established
independence and equalities \eqref{esc3} say that
$*_k\As_k(C)=\As_{n-1-k}(C^*)$ for all $k=1,...,n-2$. \proofend

\smallskip

{\bf Pl\"ucker formulas.} The degrees of the associated curves $d_k=\deg C_k$, $1\le k\le n-1$,
can be defined by the Pl\"ucker embedding $Gr(k+1,\C^{n+1})\to\Lambda^{k+1}\C^{n+1}$ or by intersection
with a generic Schubert cell $W_1=W_1(L_{n-k-1}\subset Gr(k+1,\C^{n+1})$ formed by
the $(k+1)$-subspaces intersecting with a given $(n-k-1)$-subspace
$L_{n-k-1}$ (cf. \cite[Page 268]{GH}).
They
can be computed from $d=\deg C$, $g$, and
the ramification indices $\bet_k=\bet_k(C)$, $k=0,...,n-1$,
via the Pl\"ucker formulas (see, for example, \cite[Page 270]{GH} or \cite[\S3]{P}):
\begin{equation}
d_{k-1}-2d_k+d_{k+1}=2g-2-\bet_k,\quad
k=0,...,n-1\ ,
\label{esc14}\end{equation}
where $d_{-1}=d_n=0$, $d_0=d=\deg C$.
In particular,
$$d=n(1-g)+\frac{1}{n+1}\sum_{i=0}^{n-1}(n-i)\beta_i,
\quad d^*=n(1-g)+\frac{1}{n+1}\sum_{i=0}^{n-1}(n-i)\beta_{n-1-i}\ ,$$
$$d+d^*=\sum_{i=0}^{n-1}\bet_i-2n(g-1),\quad \text{where}\quad d^*=\deg C^*\ .$$

\smallskip

{\bf Associated maps.} Let $\overline{d}=(d_0,...,d_{n-1})\in\Z^n_{>0}$. Introduce the space
\begin{equation}
{\mathcal M}_{g,0}^\nd(\PP^n,\overline{d})=\{C\in{\mathcal M}_{g,0}(\PP^n,d_0)\ :
\ \deg C_k=d_k,\ k=1,...,n-1\}\ .\label{ecf1}\end{equation}

\begin{theorem}\label{2prim}
Let $\overline d\in\Z^n_{>0}$ and ${\mathcal M}_{g,0}^\nd(\PP^n,\overline d)\ne\emptyset$.
Then there are well-defined injective maps
$$\As_k:{\mathcal M}^\nd_{g,0}(\PP^n,\overline d)\to{\mathcal M}_{g,0}(Gr(k+1,\C^{n+1}),d_k),\quad
k=1,...,n-1\ ,$$ which are birational onto their images.
Furthermore, $$\As_{n-1}:{\mathcal M}^\nd_{g,0}(\PP^n,\overline d)\to
{\mathcal M}^\nd_{g,0}((\PP^n)^*,\overline d^*),\quad \overline d^*=(d_{n-1},...,d_0)\ ,$$ is an birational homeomorphism, satisfying
$\As_{n-1}^2=\Id$.
\end{theorem}

{\bf Proof.} The statement on $\As_{n-1}$
follows from Theorem \ref{tsc}.

Suppose that $1\le k<n-1$. In fact, we need to show only the injectivity of
$\As_k$.
Let $C=[f:S\to\PP^n]$ and $C'=[f':S'\to\PP^n]$ be two non-degenerate curves
such that $\As_k(C)=\As_k(C')$. It is enough to show that if some germs $f:(S,s)
\to\PP^n$, $f':(S',s')\to\PP^n$ admit parameterizations $X(t)$, $Y(t)$, $t\in(\C,0)$, such that
\begin{equation}Y(t)\wedge Y'(t)\wedge...\wedge Y^{(k)}(t)=\mu(t)\cdot X(t)\wedge X'(t)\wedge
...\wedge X^{(k)}(t),\quad t\in(\C,0)\ ,\label{esc5}\end{equation} with
some non-vanishing scalar function $\mu(t)$, then
$Y(t)=\lambda(t)X(t)$, where $\lambda(t)$ is a scalar function. Without loss of generality,
we can assume that $X(t)$ is given by \eqref{esc4} with $\bet_0=...=\bet_{k-1}=0$, in particular,
$X\wedge X'\wedge...\wedge X^{(k)}$ does not vanish
for $t\in(\C,0)$. Relation \eqref{esc5} yields that
$Y,Y',...,Y^{(k)}\in\Span\{X,X',...,X^{(k)}\}$. Particularly,
$$Y(t)=\lambda(t)X(t)+\lambda_1(t)X'(t)+...+\lambda_k(t)X^{(k)}(t)\ .$$ It follows that
$Y'(t)$ is a linear combination of $X,...,X^{(k)}$ plus $\lambda_k(t)X^{(k+1)}(t)$.
However, since $C$ is non-degenerate and $k<n-1$, the vector $X^{(k+1)}$ is
linearly independent of $X,X',...,X^{(k)}$, and
hence $\lambda_k(t)\equiv0$. Considering subsequently $Y''(t),...,Y^{(k)}$, in the same manner
we deduce that $\lambda_1(t)\equiv...\equiv\lambda_{k-1}(t)\equiv0$. Hence, $Y=\lambda X$.
\proofend

In fact, the maps inverse to $\As_k$, $1\le k\le n-2$, can be explicitly exhibited as follows.
Let $f:S\to\PP^n$ be a non-degenerate curve, $f_k:S\to Gr(k+1,\C^{n+1})$ its $k$-th associated curve,
$1\le k\le n-1$. Let $s\in S$ and $\Lambda=f_k(s)\in Gr(k+1,\C^{n+1})$. Given a
function germ $F\in{\mathcal O}_{Gr(k+1,\C^{n+1}),\Lambda}$ vanishing at $\Lambda$, we denote by
$\ord_s(F\circ f_k)$ the order of zero of $F\circ f_k$ at $s$. For any proper linear subspace
$L\subset T_{\Lambda}Gr(k+1,\C^{n+1})$, we define the
{\bf intersection multiplicity} $(f_k\big|_{S,s}\cdot L)$ as
$$\min_F\left\{\ord_s(F\circ f_k)\ :\ F(\Lambda)=0,\ D_{\Lambda}F\big|_L=0\right\}\ .$$
Recall also that there exists a canonical identification (see, e.g.
\cite[Lecture 16]{H})
\begin{equation}T_\Lambda Gr(k+1,\C^{n+1})\simeq\Hom(\Lambda,\C^{n+1}/\Lambda)\ .\label{esc16}\end{equation}
Denote by $\pi_\Lambda:\C^{n+1}\to\C^{n+1}/\Lambda$ the standard projection.

\begin{lemma}\label{lsc3}
Fix $m=1,...,n-k-1$.

(1) There exists a unique linear subspace $L\subset T_{\Lambda}Gr(k+1,\C^{n+1})$ of dimension
$m$ for which the intersection multiplicity $(f_k\big|_{S,s}\cdot L)$ attains its maximal value. We
call this subspace the {\bf osculating $m$-space} to the curve $f_k:S\to Gr(k+1,\C^{n+1})$ at the point $s$,
and denote it by $L_m(f_k,s)$.

(2) Via the identification \eqref{esc16}
the following holds: the projective subspace
$$\PP\left(\pi_\Lambda^{-1}\Span\left(\bigcup_{\Lambda'\in L_m(f_k,s)}\Ima\Lambda'\right)\right)\subset\PP^n$$
coincides with the osculating $(k+m)$-plane to the curve $f:S\to\PP^n$ at the point $s$.
\end{lemma}

{\bf Proof.}
Consider the canonical parametrization \eqref{esc4} of the germ \mbox{$f:(S,s)\to\PP^n$}. The corresponding
parametrization $X(t)\wedge X'(t)\wedge...\wedge X^{(k)}(t)$ of the germ \mbox{$f_k:(S,s)\to Gr(k+1,
\C^{n+1})$} via elementary transformations converts a family of
$(k+1)\times(n+1)$ matrices $F_k(t)=(I_{k+1},A(t))$, where $I_{k+1}$ is the unit matrix of
size $(k+1)$, and $A=\left(a_{ij}(t)\right)_{i=0,...,k,\ j=1,...,n-k}$ is a
$(k+1)\times(n-k)$ matrix, in which
$$a_{ij}=t^{\gamma(i,j)}(a_{ij}^{(0)}+O(t)),\quad \gamma(i,j)=k+j-i+\sum_{i\le r\le k+j-1}\beta_r\ ,$$
and for all $j=1,...,n-k$,
$$a_{kj}^{(0)}=\prod_{r=1}^k\left(k+2-r+\sum_{l=r}^{k+1}\beta_l\right)\cdot
\left(\prod_{r=1}^k\left(k+1-r+\sum_{l=r}^k\beta_l\right)\right)^{-1}\ne0\ .$$
Note that the coordinate system $(a_{ij},\ i=0,...,k,\ j=1,...,n-k)$ in a neighborhood of $\Lambda$
in $Gr(k+1,\C^{n+1})$ comes from the projection
$(O,A)\mapsto A$ onto
\mbox{$T_\Lambda Gr(k+1,\C^{n+1})\simeq\Hom(\Lambda,\C^{n+1}/\Lambda)$} (cf. \cite[Example 16.1]{H}).

We have
$$F_k(t)=(I_{k+1},O)+t^{1+\beta_k}\cdot(O,A^{(0)}_1)+...+t^{n-k+\beta_k+...+\beta_{n-1}}
\cdot(O,A^{(0)}_{n-k})+O(t^\gamma)\ ,$$ where $\gamma>n-k+\beta_k+...+\beta_{n-1}$, and
$A^{(0)}_r$ contains the only non-zero diagonal $(a^{(0)}_{ij})_{i-j=k-r}$ with at least one non-zero
element $a_{kr}^{(0)}$. It follows that the maximal
value of $(f_k\big|_{S,s}\cdot L)$ over $m$-dimensional subspaces of
$T_\Lambda Gr(k+1,\C^{n+1})$ equals \mbox{$m+1+\beta_k+...+\beta_{k+m}$} and is attained
only at the space $L$ spanned by the (linearly independent) elements
$A_1^{(0)}$, ..., $A_m^{(0)}$. The space generated by the rows of these matrices lifts to
$\C^{n+1}$ as the span of the first $k+m+1$ unit vectors
$e_0,...,e_{k+m}$, and its projectivization is just the
$(k+m)$-th osculating plane of the curve $f:S\to\PP^n$ at $s\in S$.
\proofend

As a direct consequence of Lemma \ref{lsc3} one obtains

\begin{corollary}\label{csc1} Let $n\ge3$.
For any $k=1,...,n-2$ and $m=1,...,n-k-1$, there exists a
birational map
$$\As^k_m:\As_k\left({\mathcal M}_{g,0}^{\nd}(\PP^n,\overline d)\right)\dashrightarrow
\As_{k+m}\left({\mathcal M}_{g,0}^{\nd}(\PP^n,\overline d)\right)\ ,$$ which is
well-defined on a Zariski open, dense subset of $\As_k\left({\mathcal M}_{g,0}^{\nd}(\PP^n,\overline d)\right)$
and which extends
to a homeomorphism
$$\widehat\As^k_m:\As_k\left({\mathcal M}_{g,0}^{\nd}(\PP^n,\overline d)\right)\to
\As_{k+m}\left({\mathcal M}_{g,0}^{\nd}(\PP^n,\overline d)\right)\ .$$
Furthermore,
the relations $\widehat\As^k_m\circ\As_k=\As_{k+m}$ and $(\As_k)^{-1}=*\widehat\As^k_{n-k-1}$
are
satisfied.
\end{corollary}

{\bf Proof.}
Let $f:S\to\PP^n$ represent a smooth point of a component of
${\mathcal M}_{g,0}^{\nd}(\PP^n,\overline d)$, and let $s\in S$ be a generic point. Then
$\beta_i(f,s)=0$, $i=0,...,n-1$, and the canonical parametrization \eqref{esc4} takes form
$$x_0=1,\ x_1=t(1+O(t)),\ ...\ x_n=t^n(1+O(t))\ .$$ Then, in
a smooth neighborhood ${\mathcal U}\subset
{\mathcal M}_{g,0}^{\nd}(\PP^n,\overline d)$ of $f:S\to\PP^n$ we have
a parametrization $u\in{\mathcal U}\mapsto
[f_u:S_u\to\PP^n]$, a smooth section $u\in{\mathcal U}\mapsto
s(u)\in S_u$, and a universal parametrization
$t\in(\C,0)\to(S_u,s(u))$ given by
$$x_0(t)=1,\ x_1(t)=t(\xi(u)+O(t)),\ ...\ ,\ x_n=t^n(\xi_n(u)+O(t))\ ,$$ with
$\xi_1,...,\xi_n$ smoothly depending on $u\in{\mathcal U}$.
Finally, the computations in the proof of Lemma \ref{lsc3} provide a smooth parametrization of
a neighborhood of $f_k:S\to Gr(k+1,\C^{n+1})$ in
$\As_k\left({\mathcal M}_{g,0}^{\nd}(\PP^n,\overline d)\right)$ and define
the map $\As^k_m$ as a smooth analytic map in that neighborhood.
Combining this with Theorem \ref{2prim}, we obtain all statements.
\proofend

\begin{remark}\label{rsc1}
(1) It follows from Corollary \ref{csc1} that two curves $C'\subset\PP^n$ and $C''\subset(\PP^n)^*$ are dual
if they are $k$-dual for at least one $k=1,...,n-1$ (see Definition \ref{maindef}).

(2) A birational homeomorphism may not be an isomorphism, what one can see in
an elementary example of such a map between the line and a plane cuspidal cubic.
To prove (or disprove) that the maps in Corollary \ref{csc1} are isomorphisms,
one should study the behavior of these maps at singularities of
the spaces $\As_k\left({\mathcal M}_{g,0}^{\nd}(\PP^n,\overline d)\right)$.
\end{remark}

\smallskip

{\bf Integrable curves in Grassmannians.} The curves in Grassmannians,
which are associated to non-degenerate projective curves, are called {\bf integrable curves}.
They admit the following characterization.

Given a linear subspace $V\subset\C^{n+1}$ such that $\dim V>k+1$, we
obtain a natural embedding $Gr(k+1,V)\hookrightarrow Gr(k+1,\C^{n+1})$, whose image
is called a {\it Grassmannian subvariety} of
$Gr(k+1,\C^{n+1})$. A curve \mbox{$\varphi:S\to Gr(k+1,\C^{n+1})$} will be
called $(k,n)${\bf-non-degenerate}\footnote{In
this sense, a non-degenerate curve $\varphi:S\to\PP^n$ is $(0,n)$-non-degenerate.} if its image
is not contained in any proper Grassmannian subvariety of $Gr(k+1,\C^{n+1})$.

\begin{remark}\label{rnew}
Let us illustrate the difference between the non-degeneracy of curves in the projective space and the
$(k,n)$-non-degeneracy of curves in Grassmannians. Consider, for example, a cubic normal curve
$C=[f:\PP^1\to\PP^3]$. Its first associated curve $f_1:\PP^1\to Gr(2,\C^4)$ is $(1,3)$-non-degenerate
(cf. Theorem \ref{tsc2} below), but it becomes degenerate after the Pl\"ucker embedding
$Gr(2,\C^4)\hookrightarrow\PP^5$. Indeed, choosing a parametrization
$x_0=1$, $x_1=t$, $x_2=t^2$, $x_3=t^3$ for $C$, we obtain the following Pl\"ucker coordinates of
its first associated curve embedded into $\PP^5$: $\pi_{03}=3t^2$ and $\pi_{12}=t^2$, which means that
the spoken curve lies in the hyperplane $\pi_{03}-3\pi_{12}=0$.
\end{remark}

Let $t\in(\C,0)\mapsto L(t)=\Span\{X_0(t),...,X_k(t)\}\subset\C^{n+1}$ be a
parametrization of
a curve $\varphi:S\to Gr(k+1,\C^{n+1})$ in a neighborhood of a point $s\in S$. Denote
$L'(t)=\Span\{X'_0(t),...,X'_k(t)\}$ and let $\hat\varphi(s)=L(0)+L'(0)$. Clearly, the
subspace $\hat\varphi(s)\subset\C^{n+1}$ does not depend on the choice of a local parameter $t$
and a moving frame $(X_0(t),...,X_k(t))$.

\begin{theorem}\label{tsc2}
Let $n\ge3$ and $1\le k\le n-2$. A curve $\varphi:S\to Gr(k+1,\C^{n+1})$ is integrable if and only if
it is $(k,n)$-non-degenerate, and $\dim\hat\varphi(s)=k+2$ for almost all points $s\in S$.
\end{theorem}

{\bf Proof.}
If $f:S\to\PP^n$ is a non-degenerate curve, then its associated curve $f_k:S\to Gr(k+1,\C^{n+1})$ is
non-degenerate too. Furthermore, in  parametrization \eqref{esc4}, we have
$$\dim\hat f_k(s)=\dim\Span\{X(0),X'(0),...,X^{(k+1)}(0)\}=k+2\ ,$$
which holds for all points $s\in S$, where $\beta_0(f,s)=...=\beta_k(f,s)=0$, i.e., for
almost all points of $S$.

Now, let $\varphi:S\to Gr(k+1,\C^{n+1})$ be a non-degenerate curve with $\dim\hat\varphi(s)=k+2$ for
almost all points $s\in S$. We can choose a local parameter $t$ and a moving frame
$(X_0(t),...,X_k(t))$ spanning $\varphi(s_t)$ in a neighborhood of a regular point
$s_0\in S$ so that
$$\hat\varphi(s_t)=\Span\{X_0(t),...,X_k(t),X'_0(t)\} ,$$ and correspondingly
$$X'_1(t),...,X'_k(t)\in\Span\{X_0(t),...,X_k(t),X'_0(t)\}\ .$$
Observe that
\begin{equation}X''(t)\not\in\hat\varphi(s_t)
\quad\text{for almost all values}\quad t\in(\C,0)\ .\label{esc1209}\end{equation} Indeed, otherwise,
we would immediately get that $X^{(i)}_j(0)\in\hat\varphi(s_0)$ for all $j=0,...,k$ and
$i\ge0$, and hence $\varphi(s)\subset\hat\varphi(s_0)\subsetneq\C^{n+1}$ for all
$s\in S$, contrary to the non-degeneracy assumption.

We claim that the curve $\varphi:S\to Gr(k+1,\C^{n+1})$ can be uniquely restored from
the curve $\hat\varphi:S\to Gr(k+2,\C^{n+1})$. Indeed, let $\psi:S\to Gr(k+1,\C^{n+1})$ be another curve
such that $\hat\psi=\hat\varphi$, and let $\psi(s_t)=\Span\{Y_0(t),...,Y_k(t)\}$ in a neighborhood of
$s\in S$. Thus,
$$Y_i(t)=a_i(t)X'_0(t)+b_{i0}X_0(t)+...+b_{ik}X_k(t),\quad i=0,...,k\ .$$
It follows that
$$Y'_i(t)\in a_i(t)X''_0(t)+\hat\varphi(s_t)=a_i(t)X''_0(t)+
\hat\psi(s_t),\quad i=0,...,k\ ,$$ which in view of \eqref{esc1209} yields $a_0(t)\equiv...\equiv a_k(t)\equiv0$,
and hence the required claim.

The statement of Theorem follows now from the above claim and Corollary \ref{csc1} by descending
induction starting at $k=n-2$.
\proofend

\smallskip

{\bf Open problems.} It is natural to call the families ${\mathcal M}_{g,0}^{\nd}(\PP^n,\overline d)$
defined by \eqref{ecf1} the {\bf equiclassical families} of curves in $\PP^n$, $n\ge 2$.
Similarly to the planar case, we formulate

\begin{question}\label{q1}
(1) Under what conditions is ${\mathcal M}_{g,0}^{\nd}(\PP^n,\overline d)$ equidimensional
and what is $\dim{\mathcal M}_{g,0}^{\nd}(\PP^n,\overline d)$ ? Under what conditions
is ${\mathcal M}_{g,0}^{\nd}(\PP^n,\overline d)$ irreducible ? What are generic members of
(the components of) ${\mathcal M}_{g,0}^{\nd}(\PP^n,\overline d)$ ?

(2) Is the duality map $\As_{n-1}:{\mathcal M}_{g,0}^{\nd}(\PP^n,\overline d)
\to{\mathcal M}_{g,0}^{\nd}(\PP^n,\overline d^*)$, where
$\overline d^*=(d_{n-1},...,d_0)$, an isomorphism ?
\end{question}

G. Canuto \cite{Ca} has considered a problem related to the question on a generic member of
${\mathcal M}_{g,0}^{\nd}(\PP^n,\overline d)$ and proved that any compact
Riemann surface, having a non-special linear system $g^n_d$, $n\ge3$, without base points, can be mapped
to $\PP^n$ as a non-degenerate curve of degree $d$ so that its first $n-2$ associated curves are
smooth.

Following the treatment of curves in Grassmannians in \cite{Ca,Per}, one can state problems similar to
those in Question \ref{q1}. However, families of integrable curves
$${\mathcal M}_{g,0}^{\nd}(Gr(k+1,\C^{n+1}),\overline d):=
\As_k({\mathcal M}_{g.0}^{\nd}(\PP^n,\overline d)),\quad k=1,...,n-2\ ,$$
seem to be most interesting:

\begin{question}\label{q2}
(1) Are the maps $\As_k:{\mathcal M}_{g.0}^{\nd}(\PP^n,\overline d)\to
{\mathcal M}_{g,0}^{\nd}(Gr(k+1,\C^{n+1}),\overline d)$, $k=1,...,n-2$, isomorphisms ?
Are the maps $\widehat\As^k_m:{\mathcal M}_{g,0}^{\nd}(Gr(k+1,\C^{n+1}),\overline d)
\to{\mathcal M}_{g,0}^{\nd}(Gr(k+m+1,\C^{n+1}),\overline d)$, $1\le k\le n-2$,
$1\le m\le n-k-1$, isomorphisms ?

(2) What are generic members of (the components of)
${\mathcal M}_{g,0}^{\nd}(Gr(k+1,\C^{n+1}),\overline d)$, $k=1,...,n-2$ ?
\end{question}

{\ncsc Steklov Mathematical Institute \\[-21pt]

Gubkina str., 8, \\[-21pt]

Moscow 119991, Russia. \\[-21pt]

{\it E-mail address}: {\ntt kulikov@mi.ras.ru}

\vskip10pt

{\ncsc School of Mathematical Sciences \\[-21pt]

Raymond and Beverly Sackler Faculty of Exact Sciences\\[-21pt]

Tel Aviv University \\[-21pt]

Ramat Aviv, Tel Aviv 69978, Israel} \\[-21pt]

{\it E-mail address}: {\ntt shustin@post.tau.ac.il}

\end{document}